\documentclass[11pt]{article}
\usepackage[a4paper,margin=25mm]{geometry}
\usepackage[T1]{fontenc}
\usepackage[utf8]{inputenc}
\usepackage{lmodern}
\usepackage{amsmath,amssymb,mathtools,bm}
\usepackage{microtype,booktabs,array}
\usepackage{graphicx}
\usepackage{subcaption}
\usepackage{float}
\usepackage{xcolor}
\usepackage{listings}
\usepackage{enumitem}
\usepackage{hyperref}
\usepackage[nameinlink,noabbrev]{cleveref}
\usepackage[numbers,sort&compress]{natbib}
\hypersetup{colorlinks=true,linkcolor=blue!45!black,citecolor=blue!45!black,urlcolor=blue!45!black}
\newcommand{\R}{\mathbb R}
\newcommand{\Lie}{\operatorname{Lie}}
\newcommand{\SO}{\operatorname{SO}}
\newcommand{\Reach}{\mathcal R}
\newcommand{\g}{\mathfrak g}
\newcommand{\salg}{\mathfrak s}
\newcommand{\ad}{\operatorname{ad}}
\newcommand{\Span}{\operatorname{span}}
\definecolor{codegray}{gray}{0.96}
\lstdefinestyle{mma}{basicstyle=\ttfamily\footnotesize,backgroundcolor=\color{codegray},frame=single,breaklines=true,columns=fullflexible,showstringspaces=false}
\title{Nonlinear Controllability and the Propagation of Local Information\\\large From the Kalman Family to Lie Brackets, Rotation Groups, and Reachable Subgroups}
\author{Philippe Mullhaupt\\Laboratoire d'Automatique}
\date{}
\begin{document}
\maketitle

\begin{abstract}
This article develops a self-contained Lie-theoretic route from linear controllability to nonlinear controllability on matrix Lie groups. The organizing question is whether local algebraic information can be propagated into statements about reachable sets. The linear case supplies the model: the matrix exponential, Cayley--Hamilton theorem, trajectory formula and Kalman family $B,AB,\ldots,A^{n-1}B$ reduce controllability to finite-dimensional linear algebra. In the nonlinear setting, Lie brackets replace matrix powers, but the associated Lie algebra may be infinite dimensional and the finite-dimensional Lie correspondence can fail. Particular emphasis is placed on rotation groups. Skew-symmetric matrices, their exponentials, commutators, one-parameter subgroups, and the Baker--Campbell--Hausdorff formula provide a concrete laboratory for seeing how infinitesimal directions propagate on $\SO(n)$. A numerical $\SO(3)$ simulation illustrates the geometry of controlled rotations, while explicit $\SO(4)$ and $\SO(7)$ commutator sequences make the propagation of a localized control direction visible. The article then develops the attainable-subgroup argument for right-invariant systems, gives a detailed proof sequence behind the Yamabe step, and works through solvable, nilpotent and ideal examples for upper-triangular matrix algebras. The finite-dimensional mechanism is finally contrasted with Sussmann's local controllability theory for general nonlinear systems.
\end{abstract}

\tableofcontents

\section{Introduction: propagating local information}
A useful way to compare linear and nonlinear controllability is to ask how much information near a point can be propagated into a statement about finite-time reachability. In a linear time-invariant system, a finite family of directions --- $B,AB,\ldots,A^{n-1}B$ --- completely determines controllability. The deeper question is whether an analogous finite collection of local directions can determine reachable sets in a nonlinear system.

The natural nonlinear replacement of the Kalman family is the family of iterated Lie brackets of the vector fields defining the system. This replacement is exact in the linear case, but nonlinear vector fields introduce two essential difficulties: there is no nonlinear analogue of Cayley--Hamilton that forces the bracket family to terminate, and exact trajectories involve increasingly complicated iterated integral structures. The algebra generated by the vector fields may therefore be infinite dimensional.

Rotation matrices isolate the finite-dimensional mechanism in a setting where every object can be calculated explicitly. The Lie algebra $\mathfrak{so}(n)$ consists of skew-symmetric matrices, the corresponding Lie group is $\SO(n)$, the exponential map is the matrix exponential, and commutators of skew-symmetric matrices are again skew symmetric. This makes the rotation group an ideal laboratory for the slogan ``local information propagates.'' The examples below are therefore treated as part of the main argument rather than as secondary illustrations.

\section{The linear prototype}
Consider the single-input LTI system
\begin{equation}
\dot x=Ax+Bu,\qquad x\in\R^n,\quad u\in\R.
\label{eq:lti}
\end{equation}
Let $\Phi(x_0,t,u)$ denote the solution generated by an input $u(\cdot)$. The unrestricted-time reachable set from the origin is
\[
\Reach(0)=\{x\in\R^n:\exists t>0,\exists u(\cdot),\ \Phi(0,t,u)=x\},
\]
and the reachable set at a specified time is
\[
\Reach(x_0,t)=\{x\in\R^n:\exists u(\cdot),\ \Phi(x_0,t,u)=x\}.
\]
Small-time local controllability at the origin means that $\Reach(0,t)$ contains a neighborhood of the origin for every $t>0$. Global controllability means that every target state can be reached in some finite time.

\subsection{The matrix exponential and Cayley--Hamilton}
The first ingredient is the exponential
\[
e^A=I+A+\frac{A^2}{2!}+\frac{A^3}{3!}+\cdots.
\]
At first sight this introduces infinitely many matrix powers. Cayley--Hamilton removes that apparent infinity: every power $A^p$ is a linear combination of $I,A,\ldots,A^{n-1}$. Consequently $e^A$ itself is a linear combination of these first $n$ powers.

The second ingredient is the trajectory formula
\begin{equation}
\Phi(x_0,t,u)=e^{At}x_0+\int_0^t e^{A(t-\tau)}Bu(\tau)\,d\tau.
\label{eq:lineartraj}
\end{equation}
Together these facts give the familiar condition
\begin{equation}
\operatorname{rank}\,[B\ AB\ A^2B\ \cdots\ A^{n-1}B]=n.
\label{eq:kalman}
\end{equation}
For multiple inputs, all families $B_i,AB_i,\ldots,A^{n-1}B_i$ are included.

\subsection{The oscillator example}
For
\[
A=\begin{bmatrix}0&1\\-1&0\end{bmatrix},\qquad B=\begin{bmatrix}0\\1\end{bmatrix},
\]
one has
\[
[B\ AB]=\begin{bmatrix}0&1\\1&0\end{bmatrix},
\]
which has rank two. Thus the position and velocity of the oscillator can be assigned arbitrarily by a suitable control over a finite interval.

\section{Lie brackets as the nonlinear Kalman family}
Let $M$ be a smooth manifold. A vector field $f$ assigns a tangent vector $f(p)\in T_pM$ at each point. In coordinates, we use the Lie bracket
\begin{equation}
[f,g]=Dg\,f-Df\,g.
\label{eq:bracket}
\end{equation}
It is antisymmetric and satisfies Jacobi's identity, so the vector fields form a Lie algebra.

The linear case is recovered immediately. For $f(x)=Ax$ and the constant vector field $g(x)=B$,
\[
[f,g]=-AB,\qquad [f,[f,g]]=A^2B,\qquad [f,[f,[f,g]]]=-A^3B,\ldots
\]
Hence the Kalman family is exactly an iterated-bracket family, up to alternating signs.

For a scalar-input system
\begin{equation}
\dot x=f_0(x)+u f_1(x),
\label{eq:controlaffine}
\end{equation}
we introduce spaces $S_k(f,g)$ consisting of linear combinations of Lie monomials in $f$ and $g$ containing $g$ at most $k$ times. These spaces organize the bracket information needed in the local controllability condition.

\section{Nonlinear controllability and the HLCC conditions}
The following conditions are useful in the scalar-input setting, called the Hermes Local Controllability Conditions (HLCC):
\begin{enumerate}[label=(HLCC \arabic*)]
\item $x_0$ is a regular equilibrium for some admissible constant control $\bar u$, so $f_0(x_0)+\bar u f_1(x_0)=0$;
\item $\dim\Lie(f_0,f_1)(x_0)=\dim M$;
\item for the increasing sequence $S_k(f_0+\bar u f_1,f_1)(x_0)$, equality $S_k=S_{k+1}$ is required whenever $k$ is odd.
\end{enumerate}
Sussmann's 1983 theorem states that these conditions imply small-time local controllability \cite{sussmann1983}. The theorem is deliberately used as motivation rather than proved: its proof is much more elaborate than the linear rank test and exposes precisely the limitations of ordinary finite-dimensional Lie theory. Sussmann's later general theorem develops the nilpotent-approximation viewpoint further \cite{sussmann1987}.

\section{Rotation matrices: the central finite-dimensional example}
\label{sec:rotations}
The rotation-matrix sequence is not merely an illustration. It is the concrete finite-dimensional model through which the correspondence between a Lie algebra, its exponential map, and a Lie group becomes visible.

\subsection{Skew-symmetric matrices}
Consider the generic $3\times3$ skew-symmetric matrix
\begin{equation}
T(a,b,c)=\begin{bmatrix}
0&a&b\\
-a&0&c\\
-b&-c&0
\end{bmatrix},\qquad T^\top=-T.
\label{eq:Ttemplate}
\end{equation}
Three numerical instances are
\[
T_1=\begin{bmatrix}0&3&4\\-3&0&5\\-4&-5&0\end{bmatrix},\quad
T_2=\begin{bmatrix}0&6&13\\-6&0&78\\-13&-78&0\end{bmatrix},\quad
T_3=\begin{bmatrix}0&18&27\\-18&0&14\\-27&-14&0\end{bmatrix}.
\]
The Lie product is the matrix commutator
\[
[A,B]=AB-BA.
\]
For the first two matrices one obtains
\begin{equation}
[T_1,T_2]=\begin{bmatrix}
0&-247&204\\247&0&-15\\-204&15&0
\end{bmatrix},
\label{eq:Tcomm}
\end{equation}
which is again skew symmetric. Thus the class is closed under the bracket. Antisymmetry and Jacobi's identity follow directly from the commutator algebra.

\subsection{Exponentials and orthogonal matrices}
If $T^\top=-T$, then
\[
(e^T)^\top e^T=e^{T^\top}e^T=e^{-T}e^T=I.
\]
Moreover $\det e^T=e^{\operatorname{tr}T}=1$. Hence $e^T\in\SO(3)$. Direct numerical exponentiation confirms the orthogonality of rows and columns. The same construction is then repeated in dimension five, using an arbitrary $5\times5$ skew-symmetric matrix, to emphasize that nothing is special to dimension three.

This gives the first important propagation mechanism: the tangent directions represented by skew-symmetric matrices exponentiate into finite rotations. Studying the linear space $\mathfrak{so}(n)$ therefore generates information about the curved manifold $\SO(n)$.

\subsection{A one-parameter rotation subgroup}
Consider the explicit curve
\begin{equation}
G(t)=\begin{bmatrix}
\cos t&\sin t&0\\
-\sin t&\cos t&0\\
0&0&1
\end{bmatrix}.
\label{eq:Gt}
\end{equation}
Every $G(t)$ is orthogonal. Differentiating at $t=0$ gives
\[
\dot G(0)=\begin{bmatrix}0&1&0\\-1&0&0\\0&0&0\end{bmatrix}=:T_g.
\]
Direct evaluation verifies that $G(t)=e^{tT_g}$. This example is the simplest manifestation of a one-parameter subgroup: an infinitesimal generator $T_g$ determines an entire curve in the group.

\subsection{The group commutator and the BCH mechanism}
For two infinitesimal generators, the group commutator
\[
e^{\varepsilon T_1}e^{\varepsilon T_2}e^{-\varepsilon T_1}e^{-\varepsilon T_2}
\]
produces, to lowest nontrivial order, the Lie bracket direction $[T_1,T_2]$. Successive truncations of the Baker--Campbell--Hausdorff expansion give
\begin{equation}
\log(e^Ae^B)=A+B+\frac12[A,B]+\frac1{12}[A,[A,B]]+\frac1{12}[B,[B,A]]+\cdots.
\label{eq:BCH}
\end{equation}
Numerically, adding the bracket terms progressively reduces the discrepancy between the exponential of the truncated Lie series and $e^Ae^B$. The conceptual conclusion is explicit: local algebraic information propagates into the Lie group. The qualification is equally important: convergence limits how far a local BCH expansion can be used directly.

\section{Why Lie's third theorem matters --- and why it can fail}
Finite-dimensional Lie theory supplies a connected simply connected Lie group for every finite-dimensional Lie algebra. In the matrix setting the correspondence is concrete: matrix commutators live in the Lie algebra and the matrix exponential maps toward the group.

This mechanism cannot simply be transferred to the full Lie algebra of smooth vector fields. The algebra of vector fields is typically infinite dimensional. The cited example of Sergeraert \cite{sergeraert1977} shows that even diffeomorphisms infinitely tangent to the identity need not embed into a one-parameter group. Thus the heuristic ``vector field $\leftrightarrow$ one-parameter subgroup'' has genuine limitations outside the finite-dimensional setting.

A direct elementary warning is provided by the vector fields With
\[
f_1=(1,x_1^2)^\top,\qquad f_2=(x_1^2,1)^\top,
\]
repeated brackets produce terms of successively higher degree, including $x_1^5,x_1^6,x_1^7$, showing concretely how an infinite family of independent vector fields can arise.

\section{The $\SO(3)$ controlled rotation simulation}
\label{sec:so3sim}
A central example collects a nine-state nonlinear system by collecting the states into a $3\times3$ matrix
\[
X=\begin{bmatrix}
x_1&x_2&x_3\\x_4&x_5&x_6\\x_7&x_8&x_9
\end{bmatrix}.
\]
With the skew-symmetric matrices
\begin{equation}
T_1=T(3,4,5),\qquad T_2=T(6,13,78),
\label{eq:T1T2sim}
\end{equation}
the system is written
\begin{equation}
\dot X=T_1X+u(t)T_2X,
\qquad u(t)=0.2\cos t,
\qquad X(0)=I_3.
\label{eq:SO3sim}
\end{equation}
This matrix representation exposes the invariant rotation-group structure.

Since $T_1+u(t)T_2$ is skew symmetric for every time, the flow preserves orthogonality. Indeed,
\[
\frac{d}{dt}(X^\top X)=X^\top(T_1+uT_2)^\top X+X^\top(T_1+uT_2)X=0.
\]
With $X(0)=I$, one therefore has $X(t)\in\SO(3)$ throughout the simulation. This identity explains the geometry seen in the three Mathematica drawings: each row of $X(t)$ remains a unit vector and hence traces a curve on the unit sphere, while the three rows remain mutually orthogonal.

\begin{figure}[H]
\centering
\begin{subfigure}[t]{0.32\textwidth}
\centering\includegraphics[width=\linewidth]{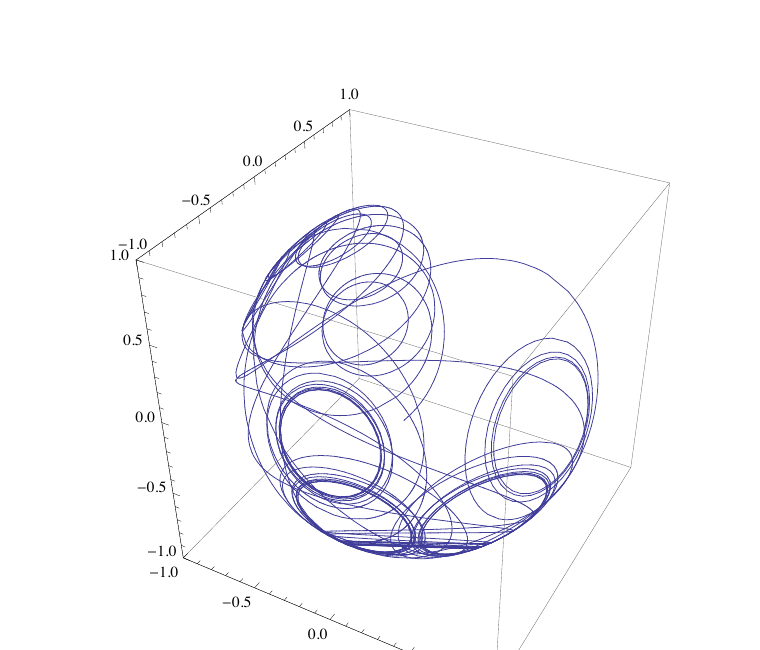}
\caption{First row of $X(t)$.}
\end{subfigure}\hfill
\begin{subfigure}[t]{0.32\textwidth}
\centering\includegraphics[width=\linewidth]{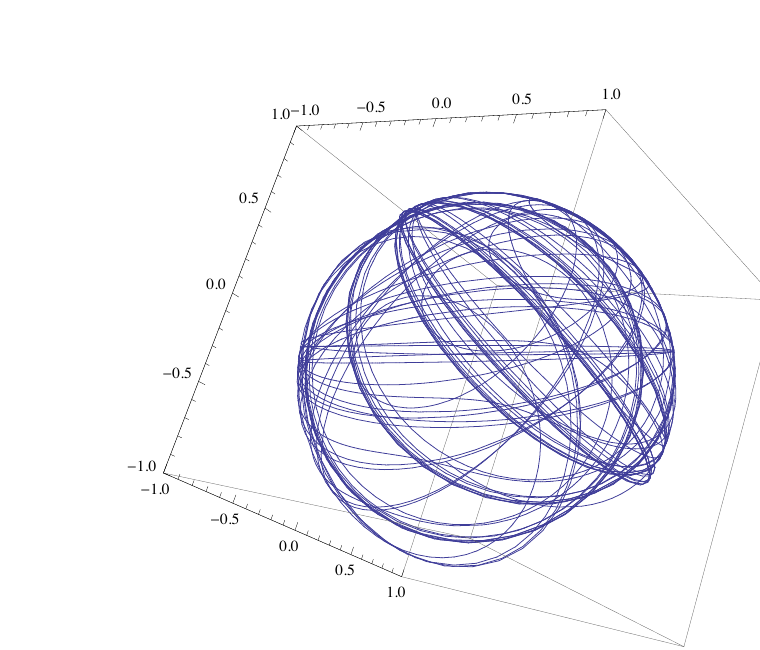}
\caption{Second row of $X(t)$.}
\end{subfigure}\hfill
\begin{subfigure}[t]{0.32\textwidth}
\centering\includegraphics[width=\linewidth]{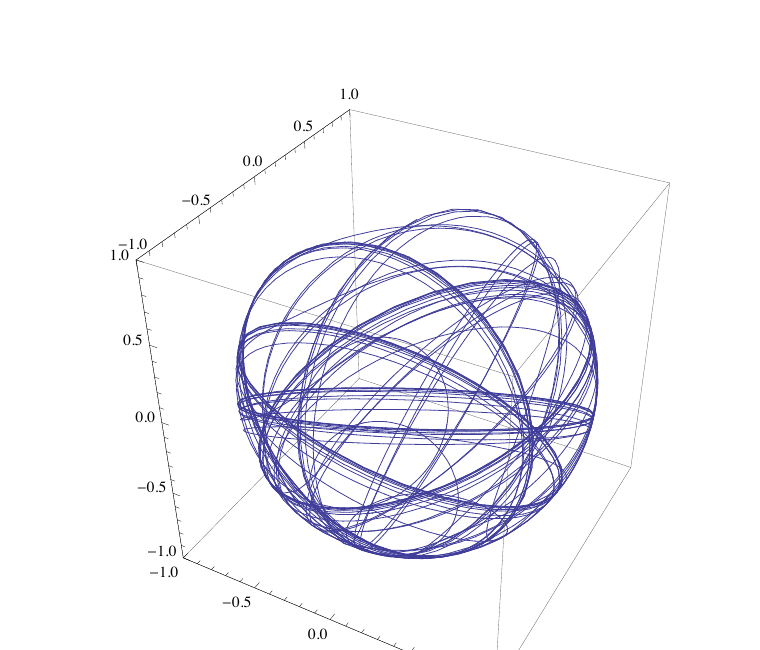}
\caption{Third row of $X(t)$.}
\end{subfigure}
\caption{The three trajectory drawings from the Mathematica simulation of \cref{eq:SO3sim}. Each row evolves on the unit sphere because the state matrix remains in $\SO(3)$.}
\label{fig:rowtrajectories}
\end{figure}

The corresponding Mathematica construction is included because it makes the simulation reproducible:
\begin{lstlisting}[style=mma,caption={Core Mathematica simulation.}]
T  = {{0,a,b},{-a,0,c},{-b,-c,0}};
T1 = T /. {a->3,b->4,c->5};
T2 = T /. {a->6,b->13,c->78};
XXt = {{x1[t],x2[t],x3[t]},
       {x4[t],x5[t],x6[t]},
       {x7[t],x8[t],x9[t]}};
xx0 = Thread[Flatten[XXt /. t->0] == Flatten[IdentityMatrix[3]]];
equDiff = Thread[Flatten[D[XXt,t]] ==
                 Flatten[(T1 + u T2).XXt]] /. u->0.2 Cos[t];
sols = NDSolve[Join[equDiff,xx0],Flatten[XXt],{t,0,30}][[1]];
ParametricPlot3D[Part[XXt,1]/.sols,{t,0,30},PlotRange->{{-1,1},{-1,1},{-1,1}}]
\end{lstlisting}
Repeating the final command for rows two and three produces \cref{fig:rowtrajectories}.

\section{Right-invariant systems on Lie groups}
Consider the right-invariant control system
\begin{equation}
\dot x(t)=X_0(x(t))+\sum_{i=1}^m u_i(t)X_i(x(t))
\label{eq:rightinv}
\end{equation}
on a finite-dimensional Lie group $G$. Let $e$ be the identity, $\Reach(e)$ the attainable set from $e$, $\salg$ the Lie subalgebra generated by $X_0,\ldots,X_m$, and $S$ the connected subgroup associated with $\salg$.

Controllability reduces to two structural questions: when is $\Reach(e)$ a subgroup, and, once it is a subgroup, when is it the whole of $G$?

\subsection{Semigroup property}
Concatenation of controls implies that $\Reach(e)$ is a semigroup. If $g$ is reached under one control and $g'$ under a second, the concatenated control reaches the corresponding product. Right invariance is the mechanism that turns concatenation of trajectory segments into multiplication in the group.

\subsection{The homogeneous case gives inverses}
For the homogeneous system
\begin{equation}
\dot x(t)=\sum_{i=1}^m u_i(t)X_i(x(t)),
\label{eq:homogeneous}
\end{equation}
time reversal together with sign reversal of the controls generates inverse motions. Therefore $\Reach(e)$ is a subgroup rather than merely a semigroup.

\subsection{Subalgebras, ideals, solvability and nilpotency}
A subalgebra $\mathfrak h$ satisfies $[\mathfrak h,\mathfrak h]\subseteq\mathfrak h$, while an ideal $\mathfrak i$ satisfies $[\mathfrak i,\g]\subseteq\mathfrak i$. The derived series
\[
\g^{(0)}=\g,\qquad \g^{(k+1)}=[\g^{(k)},\g^{(k)}]
\]
defines solvability when it reaches zero, and the lower central series
\[
\g_0=\g,\qquad \g_{k+1}=[\g,\g_k]
\]
defines nilpotency. The distinction is important for control because both constructions measure how rapidly new bracket directions disappear, but they do so with different nesting rules.

\subsubsection{Upper-triangular matrices are solvable}
Let $\mathfrak b_3$ be the vector space of all real upper-triangular $3\times3$ matrices,
\[
U=\begin{bmatrix}a&b&c\\0&d&e\\0&0&f\end{bmatrix}.
\]
For two such matrices $U_1,U_2$, the commutator has the form
\begin{equation}
[U_1,U_2]=
\begin{bmatrix}
0&\alpha&\beta\\
0&0&\gamma\\
0&0&0
\end{bmatrix},
\label{eq:upper-first-derived}
\end{equation}
where, writing the entries of $U_i$ as $a_i,b_i,c_i,d_i,e_i,f_i$,
\[
\alpha=-a_2b_1+a_1b_2-b_2d_1+b_1d_2,
\]
\[
\gamma=-d_2e_1+d_1e_2-e_2f_1+e_1f_2,
\]
and
\[
\beta=-a_2c_1+a_1c_2-b_2e_1+b_1e_2-c_2f_1+c_1f_2.
\]
Thus the first derived algebra lies in the strictly upper-triangular algebra $\mathfrak n_3$. If $V,W\in\mathfrak n_3$, then
\[
[V,W]=\begin{bmatrix}0&0&\delta\\0&0&0\\0&0&0\end{bmatrix}
\in \Span\{E_{13}\}.
\]
Since $[E_{13},E_{13}]=0$, the next derived algebra vanishes. Hence
\[
\mathfrak b_3^{(0)}=\mathfrak b_3,\qquad
\mathfrak b_3^{(1)}\subseteq\mathfrak n_3,\qquad
\mathfrak b_3^{(2)}\subseteq\Span\{E_{13}\},\qquad
\mathfrak b_3^{(3)}=0,
\]
so $\mathfrak b_3$ is solvable.

\subsubsection{Strictly upper-triangular matrices are nilpotent}
For
\[
N=\begin{bmatrix}0&a&b\\0&0&c\\0&0&0\end{bmatrix}
\in\mathfrak n_3,
\]
one computes
\[
[N_1,N_2]=
\begin{bmatrix}
0&0&a_1c_2-a_2c_1\\
0&0&0\\
0&0&0
\end{bmatrix}.
\]
Therefore
\[
[\mathfrak n_3,\mathfrak n_3]=\Span\{E_{13}\},
\qquad
[\mathfrak n_3,E_{13}]=0.
\]
The lower central series terminates after two nontrivial steps, so $\mathfrak n_3$ is nilpotent of class two. The calculation also shows why nilpotency is stronger than solvability: the lower central series brackets every new layer again with the whole algebra.

\subsubsection{A solvable upper-triangular algebra need not be nilpotent}
Take
\[
H=\operatorname{diag}(1,0,0),\qquad E=E_{12}.
\]
Then $H,E\in\mathfrak b_3$ and
\[
[H,E]=E,
\qquad
\ad_H^k(E)=E\quad\text{for all }k\ge1.
\]
Hence the lower central series of $\mathfrak b_3$ cannot terminate at zero. Thus $\mathfrak b_3$ is solvable but not nilpotent. This compact computation expresses the same phenomenon as the repeated symbolic commutators of general upper-triangular matrices: diagonal directions can keep rescaling an off-diagonal direction indefinitely.

\subsubsection{Two concrete ideals}
Two nested subspaces make the ideal condition visible. Define
\[
\mathfrak i_1=
\left\{\begin{bmatrix}a&b&c\\0&d&e\\0&0&0\end{bmatrix}\right\},
\qquad
\mathfrak i_2=
\left\{\begin{bmatrix}a&b&c\\0&0&e\\0&0&0\end{bmatrix}\right\}.
\]
For every $I\in\mathfrak i_1$ and $U\in\mathfrak b_3$, the last diagonal entry of $[I,U]$ is zero, so $[\mathfrak i_1,\mathfrak b_3]\subseteq\mathfrak i_1$. Likewise the bracket of an element of $\mathfrak i_2$ with an arbitrary upper-triangular matrix again has the structural form of $\mathfrak i_2$. Hence
\[
\mathfrak i_2\subset\mathfrak i_1\subset\mathfrak b_3
\]
are ideals. These examples are useful later because reachable-set arguments naturally generate subalgebras and ideals associated with selected control directions.

\section{From an attainable subgroup to the generated Lie subgroup}
The central finite-dimensional statement is the following: if $\Reach(e)$ is a subgroup, then it coincides with the subgroup $S$ associated with the Lie algebra generated by the system fields. This is the crucial step in translating local algebraic generation into global reachability.

The subtle point is that an abstract subgroup of a Lie group need not automatically inherit the required embedded Lie-subgroup structure. Yamabe's theorem gives the needed bridge: an arcwise connected subgroup of a Lie group is a Lie subgroup \cite{yamabe1950}. Since controlled trajectories provide arcs from the identity, $\Reach(e)$ is arcwise connected.

Let $A=\Reach(e)$ and let $\mathfrak a$ denote its Lie algebra. The inclusion $A\subseteq S$ gives $\mathfrak a\subseteq\salg$. Conversely, take constant controls with components $a_i\in\{\pm1\}$ and define
\[
X(a)=X_0+\sum_{i=1}^m a_iX_i.
\]
The positive-time curve $\exp(tX(a))$ lies in $A$. If $A$ is a group, negative times lie there as well, so $X(a)\in\mathfrak a$. The resulting family of $X(a)$ generates the original system directions, hence $\salg\subseteq\mathfrak a$. Therefore
\begin{equation}
\Reach(e)=S.
\label{eq:ReachS}
\end{equation}
This is the classical Lie-group control mechanism developed by Jurdjevic and Sussmann \cite{jurdjevic1972}.

\section{The Yamabe mechanism: from an arcwise connected subgroup to a Lie subgroup}
\label{sec:yamabe}
Yamabe's theorem states that an arcwise connected subgroup $A$ of a finite-dimensional Lie group $G$ is a Lie subgroup \cite{yamabe1950}. For controllability this result is decisive: an attainable set may first be known only as a subgroup and as a union of controlled arcs, while the Lie-algebra argument requires a genuine Lie subgroup. The following proof sequence makes the finite-dimensional mechanism explicit.

\subsection{Shrinking neighborhoods and infinitesimal directions}
Let $e$ denote the identity of $G$, and choose a nested basis of identity neighborhoods
\[
U_1\supset U_2\supset\cdots,\qquad \bigcap_{k\ge1}U_k=\{e\}.
\]
Let $C_k$ be the arcwise connected component of $U_k\cap A$ containing $e$. Because $A$ is arcwise connected, every sufficiently small element of $A$ can be connected to $e$ through a curve in $A$, and by shrinking the neighborhood one can keep such curves arbitrarily close to $e$.

Consider sequences $a_k\in C_k$, $a_k\ne e$, with $a_k\to e$. In a coordinate chart around $e$, the directions from $e$ to $a_k$ have convergent subsequences on the unit sphere. Each limiting direction determines an infinitesimal vector $X$ in $T_eG$. Let $\mathfrak a\subset T_eG$ be the aggregate of all infinitesimal vectors obtained in this way. The central task is to prove that $\mathfrak a$ is not merely a cone of limit directions but a Lie algebra.

For $X\in\mathfrak a$, let
\[
H_X=\{\exp(tX):-1\le t\le1\}
\]
denote the corresponding local one-parameter subgroup. The construction of $\mathfrak a$ implies that, for every sufficiently small neighborhood $V$ of $e$, one can find elements of $A$ and continuous joining curves in $A$ that remain arbitrarily close to $H_XV$. In this sense the one-parameter subgroup is locally shadowed by arcs contained in $A$.

\subsection{Closure under addition}
Let $X,Y\in\mathfrak a$. For large $n$, choose elements of $A$ lying close to $\exp(X/n)$ and $\exp(Y/n)$, with joining arcs contained in a sufficiently small neighborhood of $e$. Products of these small elements remain in $A$ because $A$ is a group. The classical product limit
\begin{equation}
\left(\exp(X/n)\exp(Y/n)\right)^n\longrightarrow \exp(X+Y)
\label{eq:trotter-local}
\end{equation}
shows that the corresponding product arcs accumulate on the one-parameter subgroup generated by $X+Y$. More generally, by choosing integer exponents $r_n$ with $r_n/n\to r$, one obtains
\[
\left(\exp(X/n)\exp(Y/n)\right)^{r_n}\longrightarrow \exp(r(X+Y)).
\]
Hence $X+Y$ is again an infinitesimal direction generated by elements of $A$, and therefore $X+Y\in\mathfrak a$. Scalar closure follows similarly by reparametrizing one-parameter subgroups. Thus $\mathfrak a$ is a vector subspace of $T_eG$.

\subsection{Closure under the Lie bracket}
The bracket is recovered from the group commutator. For small $t$,
\begin{equation}
\exp(-tX)\exp(-tY)\exp(tX)\exp(tY)
=\exp\!\left(t^2[X,Y]+O(t^3)\right).
\label{eq:groupcomm-bracket}
\end{equation}
Choose $t=1/n$. Since $A$ is a group, the product of approximating elements and their inverses remains in $A$. If integers $s_n$ are chosen with $s_n/n^2\to s$, then
\begin{equation}
\left(\exp(-X/n)\exp(-Y/n)\exp(X/n)\exp(Y/n)\right)^{s_n}
\longrightarrow \exp(s[X,Y]).
\label{eq:commutator-limit}
\end{equation}
The same neighborhood-control argument used for addition then shows that $[X,Y]\in\mathfrak a$. Therefore
\[
X,Y\in\mathfrak a\quad\Longrightarrow\quad X+Y\in\mathfrak a,
\qquad [X,Y]\in\mathfrak a,
\]
so $\mathfrak a$ is a finite-dimensional Lie subalgebra of $\g=T_eG$.

\subsection{The connected subgroup generated by $\mathfrak a$}
By finite-dimensional Lie theory there exists a unique connected immersed Lie subgroup $A_0\subset G$ with Lie algebra $\mathfrak a$ \cite{chevalley1946,helgason1978,knapp2002}. Choose a basis
\[
X_1,\ldots,X_s
\]
of $\mathfrak a$ and extend it to a basis
\[
X_1,\ldots,X_s,X_{s+1},\ldots,X_r
\]
of $\g$. Near the identity, canonical coordinates of the second kind give a local representation
\begin{equation}
g=\exp(t_1X_1)\cdots\exp(t_sX_s)
  \exp(t_{s+1}X_{s+1})\cdots\exp(t_rX_r).
\label{eq:second-kind-coords}
\end{equation}
The first factor lies in $A_0$; the second factor is transverse to $A_0$.

\subsection{Why $A$ cannot have a transverse local component}
Take $a_k\in C_k\subset A$ with $a_k\to e$ and write, using \cref{eq:second-kind-coords},
\[
a_k=g_kh_k,
\qquad
 g_k\in A_0,
\]
where $h_k$ is composed only of the transverse exponentials generated by $X_{s+1},\ldots,X_r$. Suppose infinitely many $h_k$ were different from $e$. By selecting elements $f_k\in A$ arbitrarily close to $g_k$ and examining $f_k^{-1}a_k$, one would obtain a limiting direction belonging to $\mathfrak a$ because the elements lie in shrinking connected pieces of $A$. On the other hand, the leading direction of $h_k$ is transverse to $\mathfrak a$ by construction. This contradiction shows that sufficiently small elements of $A$ have no transverse component. Consequently a neighborhood of $e$ in $A$ is contained in $A_0$, and since both are groups this local inclusion propagates:
\begin{equation}
A\subseteq A_0.
\label{eq:AinA0}
\end{equation}

\subsection{Recovering a neighborhood of $A_0$ inside $A$}
The reverse inclusion uses the arcs that shadow the one-parameter subgroups. For each basis vector $X_i\in\mathfrak a$, choose a continuous curve $b_i(t)\in A$, $-1\le t\le1$, which stays as close as desired to $\exp(tX_i)$. Define
\[
Q=\left\{
\exp(t_1X_1)\cdots\exp(t_sX_s): |t_i|\le1
\right\}\subset A_0
\]
and the map
\begin{equation}
F:Q\longrightarrow A,
\qquad
F\!\left(\exp(t_1X_1)\cdots\exp(t_sX_s)\right)
=b_1(t_1)\cdots b_s(t_s).
\label{eq:F-map-yamabe}
\end{equation}
When the approximating curves are chosen sufficiently close to the exact exponentials, $F$ is a small perturbation of the canonical-coordinate parametrization of $Q$. A standard local degree/invariance-of-domain argument then implies that $F(Q)$ contains a neighborhood of $e$ in $A_0$. Since $F(Q)\subset A$, one obtains a neighborhood $W$ of $e$ in $A_0$ such that $W\subset A$.

Every connected Lie group is generated by any identity neighborhood. Because $A_0$ is connected and $A$ is a subgroup containing $W$, it follows that
\begin{equation}
A_0\subseteq A.
\label{eq:A0inA}
\end{equation}
Combining \cref{eq:AinA0,eq:A0inA} yields $A=A_0$. Thus the arcwise connected subgroup $A$ carries the Lie-subgroup structure associated with the infinitesimal algebra $\mathfrak a$.

\subsection{Why this theorem matters for controllability}
The proof clarifies the control-theoretic role of Yamabe's theorem. Reachability first produces paths and a semigroup or subgroup structure. The small products \cref{eq:trotter-local} and small commutators \cref{eq:commutator-limit} extract addition and Lie brackets from these finite motions. Finite-dimensional Lie theory then reconstructs a connected subgroup from the resulting tangent algebra, while the local inclusion argument identifies that subgroup with the attainable subgroup itself. The conclusion is exactly the propagation principle needed in the controllability argument: infinitesimal bracket generation determines the connected subgroup that can be reached.

\section{Homogeneous controllability criterion}
For a homogeneous right-invariant system, $\Reach(e)$ is already a subgroup. Combining this with \cref{eq:ReachS} gives the criterion:
\begin{quote}
A homogeneous right-invariant system is controllable on $G$ if and only if $G$ is connected and the Lie algebra generated by the controlled vector fields is all of $\g$.
\end{quote}
This is the closest finite-dimensional nonlinear analogue of the Kalman rank test in this setting.

\section{The homogeneous $\SO(4)$ example}
\begin{figure}[H]
\centering
\includegraphics[width=0.98\textwidth]{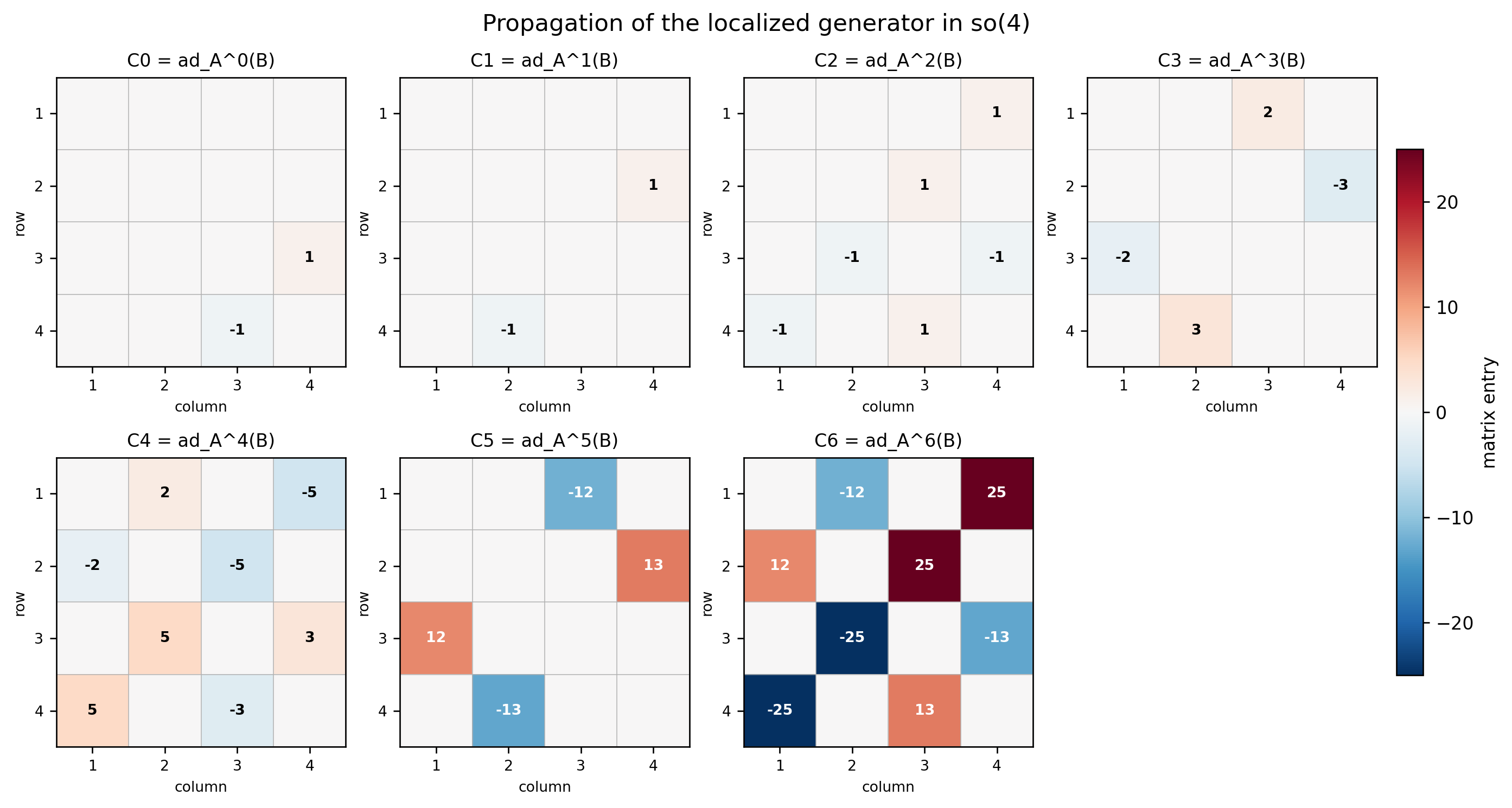}
\caption{Graphical representation of $C_k=\ad_A^k(B)$ for $k=0,\ldots,6$ in the $\SO(4)$ example. The expanding support makes the propagation of the localized rotational direction visible.}
\label{fig:so4comm}
\end{figure}

Consider
\begin{equation}
\dot X=AXu_1+BXu_2,\qquad X\in\SO(4),
\label{eq:SO4}
\end{equation}
with
\[
A=\begin{bmatrix}
0&1&0&0\\-1&0&1&0\\0&-1&0&1\\0&0&-1&0
\end{bmatrix},\qquad
B=\begin{bmatrix}
0&0&0&0\\0&0&0&0\\0&0&0&1\\0&0&-1&0
\end{bmatrix}.
\]
Both are in $\mathfrak{so}(4)$. Set
\[
C_0=B,\qquad C_{k+1}=[A,C_k]=\ad_A C_k.
\]
The first seven members are
\[
C_0=\begin{bmatrix}0&0&0&0\\0&0&0&0\\0&0&0&1\\0&0&-1&0\end{bmatrix},\quad
C_1=\begin{bmatrix}0&0&0&0\\0&0&0&1\\0&0&0&0\\0&-1&0&0\end{bmatrix},
\]
\[
C_2=\begin{bmatrix}0&0&0&1\\0&0&1&0\\0&-1&0&-1\\-1&0&1&0\end{bmatrix},\quad
C_3=\begin{bmatrix}0&0&2&0\\0&0&0&-3\\-2&0&0&0\\0&3&0&0\end{bmatrix},
\]
\[
C_4=\begin{bmatrix}0&2&0&-5\\-2&0&-5&0\\0&5&0&3\\5&0&-3&0\end{bmatrix},\quad
C_5=\begin{bmatrix}0&0&-12&0\\0&0&0&13\\12&0&0&0\\0&-13&0&0\end{bmatrix},
\]
\[
C_6=\begin{bmatrix}0&-12&0&25\\12&0&25&0\\0&-25&0&-13\\-25&0&13&0\end{bmatrix}.
\]
The algebraic pattern is clearer in graphical form. In \cref{fig:so4comm}, each cell shows a matrix entry; white cells are zero, while the signed nonzero entries are colored according to magnitude. The initial control direction $B=C_0$ is supported only on the $(3,4)$ coordinate plane. The first commutator moves this support to the $(2,4)$ plane, and subsequent commutators populate further coordinate pairs until all six independent rotational planes of $\mathfrak{so}(4)$ appear across the generated family.

This is the rotation-group analogue of the linear sequence $B,AB,A^2B,\ldots$. The relevant linear operator on the Lie algebra is now
\[
\ad_A:\mathfrak{so}(4)\to\mathfrak{so}(4),\qquad C\mapsto[A,C].
\]
The family generated from $B$ by repeated application of $\ad_A$ is therefore the natural finite-dimensional bracket counterpart of the Kalman family. In particular, the graphical sequence shows not merely growth of numerical coefficients but propagation of support through distinct coordinate planes.

\section{The $\SO(7)$ propagation example}
\begin{figure}[H]
\centering
\includegraphics[width=0.98\textwidth]{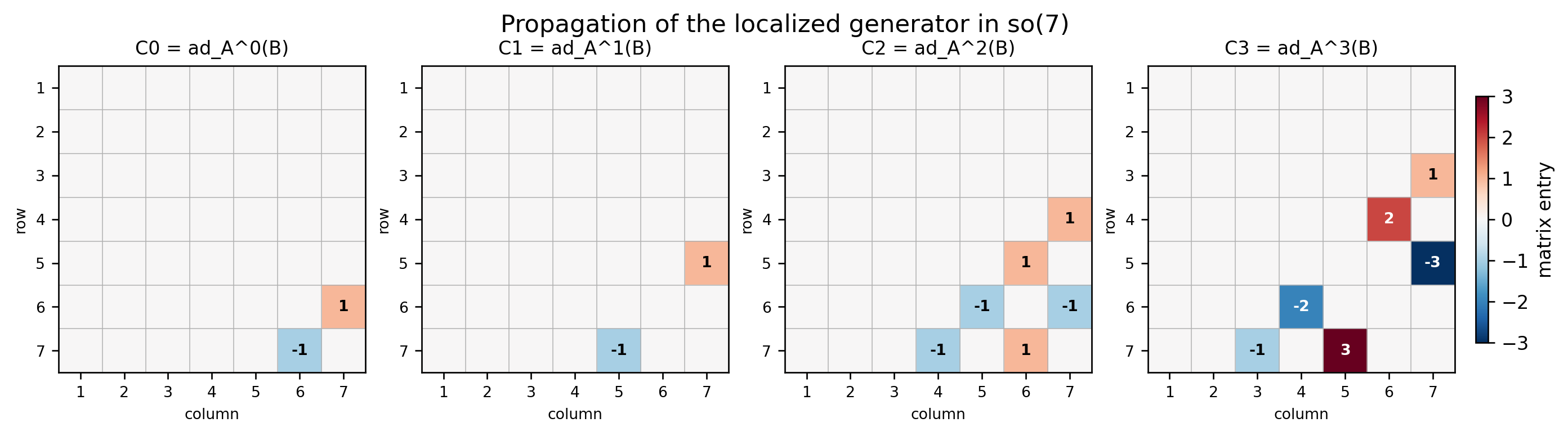}
\caption{Graphical representation of $D_k=\ad_{A_7}^k(B_7)$ for $k=0,\ldots,3$. A control generator initially confined to the $(6,7)$ plane propagates leftward through the nearest-neighbor chain under repeated commutation with $A_7$.}
\label{fig:so7comm}
\end{figure}

The same construction extends naturally to $\SO(7)$. Let $A_7$ be the skew-symmetric nearest-neighbor chain
\[
A_7=\begin{bmatrix}
0&1&0&0&0&0&0\\
-1&0&1&0&0&0&0\\
0&-1&0&1&0&0&0\\
0&0&-1&0&1&0&0\\
0&0&0&-1&0&1&0\\
0&0&0&0&-1&0&1\\
0&0&0&0&0&-1&0
\end{bmatrix},
\]
and localize the control generator in the final $(6,7)$ plane,
\[
B_7=E_{67}-E_{76}.
\]
Define $D_0=B_7$ and $D_{k+1}=[A_7,D_k]$. The first four members are
\[
D_0=E_{67}-E_{76},\qquad
D_1=E_{57}-E_{75},
\]
\[
D_2=(E_{47}-E_{74})+(E_{56}-E_{65})-(E_{67}-E_{76}),
\]
\[
D_3=(E_{37}-E_{73})+2(E_{46}-E_{64})-3(E_{57}-E_{75}).
\]
Thus the bracket depth has a direct spatial interpretation along the chain: $D_0$ acts only at the end, $D_1$ reaches one coordinate farther, $D_2$ reaches the fourth coordinate, and $D_3$ reaches the third. The full matrix values are shown graphically in \cref{fig:so7comm}.

This example makes ``propagating local information'' literal. Locality refers to the support of $B_7$ in a single rotational plane; propagation refers to the successive appearance of nonzero skew-symmetric entries in coordinate planes farther from that end. The same mechanism can be continued to higher bracket depth, and in a controllability test one asks whether the Lie algebra generated by the available directions spans all of $\mathfrak{so}(7)$.

\section{Compact groups and the nonhomogeneous case}
For compact Lie groups, compactness can restore group-like reachability properties even when a drift $X_0$ is present. In particular, in the compact Lie-group case one obtains a result of the same flavor as in the homogeneous case, whereas for a general nonhomogeneous system additional hypotheses are needed.

It is useful to distinguish three levels: finite-dimensional homogeneous systems, finite-dimensional nonhomogeneous/decomposable systems, and general infinite-dimensional vector-field algebras. The first level is completely controlled by the generated Lie algebra; the later levels require progressively more refined tools.

\section{Discussion: where the analogy succeeds and where it breaks}
The linear and finite-dimensional Lie-group theories share a common architecture. In the linear case, Cayley--Hamilton ensures that the directions generated by repeated action of $A$ form a finite family. On $\SO(n)$ and other finite-dimensional Lie groups, closure under Lie brackets and the Lie correspondence play the analogous role. BCH and group commutators explain how directions in the Lie algebra are converted into finite motions in the group.

The $\SO(3)$ simulation makes this correspondence geometric: skew-symmetric instantaneous generators preserve the orthogonal constraint, so the trajectory remains on the group and each matrix row remains on the sphere. The $\SO(4)$ and $\SO(7)$ examples make the algebraic propagation visible: a localized rotational input direction spreads through iterated commutators.

For a general nonlinear system, however, the Lie algebra of vector fields may be infinite dimensional. There is then no Cayley--Hamilton theorem to terminate the bracket generation and no unrestricted finite-dimensional Lie III correspondence to turn the algebra into a globally manageable group. This is precisely the point at which Sussmann's formal Lie series, filtrations, nilpotent approximations, and symmetry arguments become necessary \cite{sussmann1983,sussmann1987}.

\section{Conclusion}
The results above answer a common question at several levels: how can information known infinitesimally be propagated into a statement about finite-time reachability?

For LTI systems, the answer is the Kalman family and Cayley--Hamilton. For finite-dimensional matrix Lie groups, the answer is the Lie algebra, the exponential map, group commutators, BCH, and the subgroup generated by the controlled directions. The rotation examples are the most concrete expression of this mechanism: skew-symmetric matrices generate rotations; their commutators generate new infinitesimal rotations; periodic controlled dynamics remain on $\SO(3)$; and repeated adjoint operations on $\SO(4)$ and $\SO(7)$ visibly propagate a localized actuation direction through the algebra.

The analogy breaks when the generated vector-field algebra is infinite dimensional. The resulting failure is not incidental but structural, and it explains the technical depth of general nonlinear local controllability theorems. The finite-dimensional rotation examples therefore serve a double role: they provide exact controllability results in their own right and expose, in the cleanest possible setting, the mechanism that more general nonlinear theory attempts to recover by approximation.

\appendix
\section{Mathematica snippets from the rotation examples}
The following snippets record key calculations underlying the rotation-matrix examples.

\begin{lstlisting}[style=mma,caption={Skew-symmetric matrices and commutator.}]
T = {{0,a,b},{-a,0,c},{-b,-c,0}};
T1 = T /. {a->3,b->4,c->5};
T2 = T /. {a->6,b->13,c->78};
T3 = T /. {a->18,b->27,c->14};
L[A_,B_] := A.B-B.A;
L[T1,T2] // MatrixForm
\end{lstlisting}

\begin{lstlisting}[style=mma,caption={Explicit one-parameter subgroup.}]
Gg = {{Cos[t],Sin[t],0},{-Sin[t],Cos[t],0},{0,0,1}};
Tg = T /. {a->1,b->0,c->0};
D[Gg,t] /. t->0
MatrixExp[Tg 27456] - (Gg /. t->27456)
\end{lstlisting}

\begin{lstlisting}[style=mma,caption={Adjoint propagation in the $\SO(4)$ example.}]
A = {{0,1,0,0},{-1,0,1,0},{0,-1,0,1},{0,0,-1,0}};
B = {{0,0,0,0},{0,0,0,0},{0,0,0,1},{0,0,-1,0}};
MatrixForm[#]& /@ NestList[(A.#-#.A)&,B,6]
\end{lstlisting}

\end{document}